\documentclass[11pt]{amsart}
\usepackage{palatino, mathpazo, amsfonts, mathrsfs}
\usepackage{epic, eepic}
\usepackage[mathscr]{eucal}
\usepackage[all]{xy}

\newtheorem{theorem}{Theorem}[section]
\newtheorem{lemma}[theorem]{Lemma}
\newtheorem{proposition}[theorem]{Proposition}
\newtheorem{corollary}[theorem]{Corollary}

\theoremstyle{remark}

\numberwithin{equation}{section}

\newcommand{\Mbar}{\overline{\mathcal{M}}}
\newcommand{\vir}{\operatorname{vir}}

\newcommand{\Q}{\mathbb{Q}}

\newcommand{\C}{\mathbb{C}}

\newcommand{\ka}{\kappa}
\newcommand{\proj}{\mathbb P}

\newcommand{\cm}{{\mathcal{M}}}

\newcommand{\cmbar}{\overline{\cm}}

\newcommand{\al}{\alpha}

\newcommand{\be}{\beta}

\newcommand{\Ga}{\Gamma}
\newcommand{\de}{\delta}

\newcommand{\la}{\lambda}

\title{Algebraic structures on the topology of moduli spaces of curves and maps}

\author{Y.-P.~Lee} 
\address{Department of Mathematics, University of Utah, Salt Lake City, UT 84112-0090}
\email{yplee@math.utah.edu}
\author{R.~Vakil}
\address{Department of Mathematics, Stanford University, Stanford, CA~94305-2125}
\email{vakil@math.stanford.edu}
\date{September 10, 2008.}

\begin{document}

\begin{abstract}
We discuss selected topics on 
the topology of moduli spaces of curves and maps, 
emphasizing their relation with Gromov--Witten theory 
and integrable systems.
\end{abstract}

\maketitle

\setcounter{tocdepth}{1}
\tableofcontents

\setcounter{section}{-1}

\section{Introduction}

  Forty years ago, Deligne and Mumford introduced their celebrated
  compactification of the moduli space of curves.  In 1983, Mumford
  \cite{mumford} initiated a comprehensive study of intersection
  theory on the moduli space of curves, and in particular introduced
  the tautological ring.  In many ways motivated by Witten's
  conjecture \cite{witten}, the last twenty years have seen a steadily
  growing understanding of rich algebraic structures on the cohomology
  of moduli spaces of curves, and related spaces, such as
  Gromov--Witten invariants and Hurwitz numbers.  Even when
the structures deal purely with the moduli space of curves, Gromov--Witten theory
has provided a powerful tool for understanding them.
The purpose of this article is to discuss some of these structures.

{\bf 1.}  {\em Integrable systems.}
Witten's conjecture determines all top intersections of
$\psi$-classes on $\cmbar_{g,n}$, by showing that their generating
function satisfies differential equations coming from integrable
systems (the KdV hierarchy, or the Virasoro algebra).  
More precisely, let $$F_g = \sum_{n \geq 0} \frac  1 {n!} \sum_{k_1, \dots, k_n} 
\left( \int_{\cmbar_{g,n}} \psi_1^{k_1} \cdots \psi_n^{k_n} \right)
t_{k_1} \cdots t_{k_n},$$
and let $F = \sum F_g \hbar^{2g-2}$ be the generating function for all genera.
Then the KdV form of Witten's conjecture is:
\begin{eqnarray*}& & (2n+1) \frac {\partial^3} {\partial t_n \partial t_0^2} F = \\
& & \quad \left(  \frac {\partial^2} {\partial t_{n-1} \partial t_0} F \right)
\left( \frac {\partial^3} {\partial t^3_0} F \right)  +
2 \left(  \frac {\partial^3} {\partial t_{n-1} \partial t_0^2} F
\right)  \left( \frac {\partial^2 }  {\partial t_0^2} F \right)
+ \frac 1 4 \frac { \partial^5} {\partial t_{n-1} \partial t_0^4} F.
\end{eqnarray*}
 Define differential operators ($n \geq -1$)
\begin{eqnarray*}
L_{-1} &=& - \frac {\partial} {\partial t_0} + \frac {\hbar^{-2}} 2 t_0^2
+ \sum_{i=0}^{\infty}  t_{i+1} \frac \partial {\partial t_i}
\\
L_0 &=& - \frac 3 2 \frac \partial {\partial t_1} + \sum_{i=0}^{\infty}
\frac {2i+1} 2 t_i \frac {\partial} {\partial t_i} + \frac 1 {16}
\\
L_n &=&
\sum_{k=0}^{\infty} \frac { \Ga(m+n+ \frac 3 2) }{ \Ga(k+ \frac 1 2)} (t_k - \de_{k,1})
\frac \partial {\partial t_{n+k}}  \\
& & \quad \quad \quad + \frac {\hbar^2} 2 \sum_{k=1}^{n-1} (-1)^{k+1} \frac {\Ga(n-k+ \frac 1 2)}
{\Ga( -k-\frac 1 2)} \frac \partial {\partial t_k} \frac \partial {\partial t_{n-k-1}} \quad \quad (n>0)
\end{eqnarray*}
These operators satisfy
$[L_m, L_n] = (m-n) L_{m+n}$. 
The Virasoro form of Witten's conjecture, due to Dijkgraaf, Verlinde, and
Verlinde, is:  $L_n e^F = 0$ for all $n$.
These relations inductively determine the coefficients of $F$, and hence
compute all intersection numbers.

This unexpected
relationship with integrable systems has since appeared repeatedly,
and we now  have a modest understanding  as to why this relationship
might exist.  This subject is well-covered in the literature, so we content
ourselves with a short discussion in \S \ref{s:is}.

{\bf 2.}  {\em Faber-type Gorenstein conjectures.}  Based on numerical data found
using Witten's conjecture, Faber made a remarkable conjecture
\cite{faber} on the tautological part of the cohomology of
$\cm_g$. Again, this structure was quite unexpected, and in some sense
we still have little understanding ``why'' such a structure should
exist.  This conjectural structure seems to be shared by related
moduli spaces.  This topic is also well-covered (see for example
Pandharipande's ICM talk \cite{pandharipande}), so we describe
just enough  in \S \ref{s:faber} to motivate later discussion.

{\bf 3.}  {\em Polynomiality.} Proof of parts of the above conjectures
often come through another ``polynomiality'' structure of certain
invariants or intersection numbers, which is exemplified in their
first appearance, the ELSV formula \cite{elsv1, elsv2}.  For example,
the polynomiality of the ELSV formula is an easy-to-state fact
about factoring permutations into transpositions in the symmetric group,
but the only known proof requires the geometry of the moduli space
of curves.  There are other instances that look similar
that are even less understood.  In \S \ref{s:pol}, we discuss
the notion of polynomiality, and describe and explain a conjecture of Goulden, Jackson,
and the second author that is sketched in \cite{gjv}:  that 
there should be an ELSV-type formula relating 
certain double Hurwitz numbers to the intersection theory of some moduli
space.  We then describe work of Shadrin and Zvonkine \cite{sz}
and Shadrin \cite{s}
that these conjectural intersection numbers (arising from actual
Hurwitz numbers) have much richer structure than suspected in \cite{gjv}.

{\bf 4.}  \emph{Family topological field theory.} 
Madsen and Weiss' celebrated proof \cite{MW} of Mumford's conjecture 
on the stable (tautological) cohomology of the moduli space of curves has had
an unexpected and powerful application to Gromov--Witten theory.
Teleman has classified all cohomological 2D field theories based on a
semi-simple complex Frobenius algebra $A$: they are controlled by a
linear combination of $\ka$-classes and by an extension datum to the
Deligne-Mumford boundary.  This leads to a proof of Givental's
conjecture that, roughly speaking, higher genus invariants are
determined by genus $0$ invariants.  This in turn implies the Virasoro
conjecture for manifolds with semisimple quantum cohomology.
We outline this important work in some detail in \S \ref{s:teleman}.

{\bf 5.} \emph{Witten's conjecture on $r$-spin curves.}
In the last section, we briefly outline Faber--Shadrin--Zvonkine's proof 
\cite{FSZ} of Witten's conjecture on $r$-spin curves, emphasizing on the parts 
which are most related to the tautological rings.
Although this result is a consequence of Teleman's theorem, it uses 
completely different ingredients and the proof itself is quite
interesting in its own right.

Two ingredients in their proof will be briefly explained:
\begin{enumerate}
\item The first author's theorem that tautological relations hold for 
  generating functions of axiomatic (semisimple) theories \cite{ypL2}.
\item Faber--Shadrin--Zvonkine's reconstruction theorem which states that
  the Witten's correlators in higher genus can be uniquely reconstructed
  from genus zero.
\end{enumerate}
Along the way, we also give some explanation of Givental's formalism 
(\cite{aG1, LP} and references therein) which has proved important
in understanding the algebraic structure of Gromov--Witten theory.
Note that the proof of the first result can be reduced to an elementary 
statement in the geometry of boundary divisors on moduli of curves via 
Givental's theory.
See Section~\ref{s:5.3}.

{\bf Notation and background.}
We work over $\C$.
General background may be found in \cite[Ch.~23-25]{hori}, \cite{cime},
and \cite{ypL4}. 
We assume the reader is familiar with the following notions:
\newline
{\em (i)} the moduli space of pointed curves $\cm_{g,n}$ and its partial
compactifications $\cm_{g,n}^{rt}$ (curves with rational tails),
$\cm_{g,n}^{ct}$ (curves of compact type), and $\cmbar_{g,n}$ (Deligne-Mumford stable
curves, the topic of this volume).
\newline
{\em (ii)} the cohomology classes $\psi_i$ ($1 \leq i \leq n$), $\ka_j$, and
$\la_k$ ($0 \leq k \leq g$), and the tautological ring.
\newline
We also note that \cite{s2} contains an excellent review of
basic notions of cohomological field theories and Givental theory.

\subsection*{Acknowledgements}
We would like to thank T.~Jarvis, K.~Liu, S.~Shadrin, C.~Teleman, 
and D.~Zvonkine for their helpful comments on earlier versions of the paper.

Both authors are partially supported by the NSF.

\section{Integrable systems}
\label{s:is}
We content ourselves with a brief overview of
Witten's conjecture and related topics,
sufficient
to set the stage for later sections. 

It was only since the advent of Witten's conjecture that we are
able to compute all top intersections of $\psi$-classes in $\cmbar_{g,n}$, something
we now take for granted.  These allow one to compute all top intersections in the tautological ring, \cite{Falg}.  It is also important when considering Witten's conjecture
to know not just {\em that} it is true, but also {\em why} it is true, i.e.\ why
in retrospect integrable systems should be expected to control these top intersections.
Of course, Witten's original heuristic argument is one explanation.

There are now a large number of proofs of Witten's conjecture, and it
is a sign of the richness of the conjecture that almost every proof
has been fundamentally new, with fundamentally new insights.
Kontsevich's original proof (\cite{k}, see also \cite{l}) remains
unlike the others.  Mirzakhani's proof \cite{mirzakhani} (the third)
gives an elegant interpretation of each summand in the Virasoro
version of Witten's conjecture.

  The remaining proofs pass through the
ELSV formula (discussed in \S \ref{s:pol}), counting branched covers,
and fundamentally about factoring permutations in the symmetric group.
Okounkov and Pandharipande's proof \cite{op} was part of their massive
program in Gromov--Witten theory, including their proof of the
Virasoro conjecture in dimension $1$. 

 By a careful algebraic argument,  Kim and K. Liu
\cite{kimliu} give a direct proof of Witten's conjecture through
localization.      Chen, Li, and Liu thereafter gave a different
and very short proof, using ideas from \cite{lambdag} on the $\la_g$-conjecture
(a sister conjecture to Witten's conjecture, see the next section).

Meanwhile, Kazarian and Lando also gave an algebro-geometric proof of
Witten's conjecture \cite{kazarianlando}.  Shortly after that, also using
ideas from \cite{lambdag}, Kazarian gave a greatly streamlined proof \cite{kazarian}.

At each stage, our understanding of the structure behind Witten's conjectured
deepened, and in some sense we now have quite a strong understanding of why it is true.

One test of our understanding is how well it generalizes to other
situations.  We briefly mention two generalizations of Witten's
conjecture that will be relevant shortly.

First, Witten gave a
generalization of his conjecture to the moduli space of $r$-spin
curves \cite{witten2}.  Roughly speaking, an $r$-spin curve with marked
points is an $r$th root of the canonical bundle twisted by given
multiples of the points; see \S \ref{wc} for more detail.

The Virasoro conjecture (mentioned earlier), due to Eguchi, Hori,
Xiong, and also S. Katz gives Virasoro-related constraints on the
Gromov--Witten invariants of a complex projective manifold \cite[\S
10.1.4]{ck}.  The case of a point is the Virasoro form of Witten's
conjecture.

More relatives of Witten's conjecture will appear in the next section.

\section{Faber-type Gorenstein conjectures, and  three sisters}
\label{s:faber}

Faber's conjectures are well exposed in the literature; for example, the second
author has discussed them in detail in \cite[\S 3.2]{cime}.  So again
we content ourselves with saying just enough to continue our story.

\subsection{Faber's conjectures}
Faber's conjectures \cite{faber} describe an unexpected and in many
ways still-unexplained structure on the tautological ring
$R^*(\cm_g)$, which may be considered as a subring of the cohomology
ring or the Chow ring (although this technically yields two different
conjectures).

Informally speaking, Faber's conjectures on $\cm_g$ state that
$R^*(\cm_g)$ behaves like the cohomology ring of a $(g-2)$-dimensional
complex projective manifold.  Somewhat more precisely, this means:

{\bf I.}  {\em ``Vanishing/socle'' conjecture.}  $R^i(\cm_g) = 0$ for $i>g-2$, and
$R^{g-2}(\cm_g) \cong \Q$.  This was proved by Looijenga \cite{loo} and Faber
\cite[Thm.~2]{faber}.

{\bf II.}  {\em Perfect pairing conjecture.}  The analogue of
Poincar\'e duality holds: for $0\leq i \leq g-2$, the natural map
$R^i(\cm_g) \times R^{g-2-i}(\cm_g) \rightarrow R^{g-2}(\cm_g) \cong
\Q$ is a perfect pairing.  This is currently open.

{\bf III.}  {\em Intersection number conjecture.}  
Faber gives a formula for top (i.e.\ total degree $g-2$) intersections of generators
of the tautological ring, as a multiple of a generator of $R^{g-2}(\cm_g)$.
As we discuss in \S \ref{twothree}, this is proved.

The three conjectures above completely determine the tautological ring $R^*(\cm_g)$.

There are three ``sister'' conjectures that parallel Faber's, on three different spaces.
See \cite{pandharipande} for more on the three sisters, and a more detailed history.

\subsection{The first sister:  $\cmbar_{g,n}$}
On $\cmbar_{g,n}$, there is an analogous set of conjectures,  with
$g-2$ replaced by $3g-3+n$. 
This was first asked as a question by Hain and Looijenga \cite[Question~5.5]{hlo};
first stated as a speculation by Faber and Pandharipande \cite[Speculation~3]{FabP1} (in the
case $n=0$), and first stated as a conjecture by Pandharipande \cite[Conj.~1]{pandharipande}.
 In the cohomology ring, {\bf
  I} and {\bf II} trivially hold, but in the Chow ring, this is far
from clear.  ({\bf I} was first shown in \cite{gv}.)  Witten's
conjecture should be considered the intersection number part {\bf
  III}.  

\subsection{The second sister:  $\cm_g$}
\label{twothree}The pointed version of Faber's conjecture $\cm_g$ isn't for
$\cm_{g,n}$; it should be for $\cm_{g,n}^{rt}$, the space of curves
with ``rational tails'' (those stable curves with a component of
geometric genus $g$).  In this case, the role of $g-2$ is replaced by
$g-2+n$.  The intersection number portion of the conjecture (often
called {\em Faber's intersection number conjecture}) may be stated as
follows: If all $a_i > 0$, then
\begin{equation}\label{Finc}
\pi_* \psi_1^{a_1} \cdots \psi_n^{a_n} = 
\frac { (2g-3+n)! (2g-1)!!} { (2g-1)! \prod_{j=1}^n (2a_j-1)!!}
\ka_{g-2} \quad \quad \text{for $\sum a_i = g-2+n$.}
\end{equation} 
(This implies part {\bf III} of Faber's conjecture, see above, hence
we use the same name.)
This is now a theorem.  Getzler and Pandharipande showed that the
statement is a formal consequence of the Virasoro conjecture for
$\proj^2$ \cite{GeP}, in fact the large volume limit.  Givental
thereafter described how to prove the Virasoro conjecture for
projective space, and more generally Fano toric manifolds \cite{Gi1,
  aG1}.  This is a powerful approach to Faber's conjecture, but
perhaps somewhat roundabout.  Recently, K. Liu and Xu have given an
stunningly short and direct proof \cite{liuxu}.  Their approach is
quite different, and is part of their larger program for understanding
the algebraic structure of these intersection numbers. Most notably,
their approach yields a surprisingly simple explicit formula for
Witten's $n$-point function, and this has produced a number of
interesting new results about intersection numbers, of which Faber's
intersection number conjecture is just one.

As with Witten's conjecture, the susceptibility of this problem to
different approaches illustrates the richness of the question.  One
should hope that the story is not yet over, and more results should be
obtainable from the successful earlier approaches.  Furthermore, the
beautiful form of equ.\ \eqref{Finc} clearly suggests that there
should be a strong reason for it, and that further understanding
should be sought.  (The connection to double Hurwitz numbers, explored
in \cite{gjv3}, seems to be deeply related to this question.)

\subsection{The third sister: the $\la_g$-conjecture, for $\cm_{g,n}^{ct}$}
The third sister of Faber's conjecture is for the space
$\cm_{g,n}^{ct}$, the space of curves of compact type (stable curves
whose dual graphs have no loop; those stable curves with compact
Jacobian).  The corresponding conjecture is due to Faber and
Pandharipande (\cite[Spec.~2]{FabP1}, \cite[Conj.~1]{pandharipande}),
with $g-2$ is replaced by $2g-3$.  The socle portion was first proved
as a consequence of Theorem $\star$ in \cite{thmstar}.  

The intersection number portion of the conjecture was first conjecture
by Getzler and Pandharipande \cite{GeP}.  It is called the $\la_g$-conjecture because of its incarnation as a statement about
intersections on $\cmbar_{g,n}$ (rather than $\cm_{g,n}^{ct}$): For
$n,g \geq 1$,
$$
\int_{\cmbar_{g,n}} \psi_1^{b_1} \cdots \psi_n^{b_n} \la_g = \binom {2g-3+n} {b_1, \dots, b_n} c_g$$
where $\sum_{i=1}^n b_i = 2g-3+n$, $b_1, \dots, b_n \geq 0$, and $c_g$ is a 
constant depending only on $g$.

The $\la_g$-conjecture has intriguingly proved more tractable than its
sister intersection number conjectures, and a number of proofs now exist.
Faber and Pandharipande gave the first proof \cite{FPlambdag}.  A
second proof is the same as the Getzler-Pandharipande-Givental proof
of Faber's intersection number conjecture: when making the conjecture,
Getzler and Pandharipande showed that it is a formal consequence of
the Virasoro conjecture for $\proj^1$, which was later proved by
Givental (and also by Okounkov and Pandharipande).  K. Liu,
C.-C.~M.~Liu, and Zhou gave a new proof \cite{llz} as a consequence of their
proof of the Mari\~no-Vafa formula.  Finally, Goulden, Jackson, and
the second author gave a short direct (Gromov--Witten-free) proof in
\cite{lambdag} by exploiting the ``polynomiality'' structure described
in the next section, using the ELSV-formula (also described shortly).

Kazarian and, independently at the same time, Kim and K. Liu showed
that the algebraic structure introduced in \cite{gjv3}, properly
understood, also yield proofs (and explanations) of Witten's
conjecture.  (This insight was certainly not known to the
authors of \cite{gjv3}.)  Their two proofs (the most recent and
shortest proofs, mentioned above) are quite distinct, and
insightful. A complete understanding of the algebraic structures
underlying Witten's conjecture would presumably involve putting these
two proofs into a common larger framework.

\section{Polynomiality}
\label{s:pol}

We next describe the phenomenon of polynomiality of quasi-enumerative
problems on moduli spaces of curves and maps.   The central motivating
example is the ELSV formula.  

\subsection{The ELSV formula}
Fix a genus $g$, a degree $d$, and a partition of $d$ into $n$ parts,
$\al_1+ \cdots + \al_n = d$, and let $r = 2g+d+n-2$.  Fix $r+1$
distinct points $p_1, \dots, p_n, \infty$ on $\proj^1$.  Define the
{\em Hurwitz number} $H^g_{\al}$ as the number of branched covers of
$\proj^1$ by a (connected) Riemann surface, that are unbranched away
from $p_1$, \dots, $p_n$, $\infty$, such that the branching over
$\infty$ is given by $\al_1$, \dots, $\al_n$ (the monodromy lies in
the conjugacy class corresponding to that partition), and the
branching over each $p_i$ is $2+ 1+ \cdots + 1=d$ (the simplest
nontrivial branching).  We consider the $n$ preimages of $\infty$ to
be labeled.  

Up to a straightforward combinatorial factor, this corresponds to the answer
to the following combinatorial problem: given a permutation in
conjugacy class $\al$, in how many ways can it be factored into $r$
transpositions that ``connect'' the numbers $1$ through $n$ ({\em
  transitive} factorizations).  (If the condition of transitivity
seems unnatural, it is straightforward to connect to this to the
problem without the transitive condition.  This is equivalent to
counting potentially disconnected covers.  The algebraically simplest
way to relate them: the exponential of the generating function
counting connected covers is the generating function counting
potentially disconnected covers.)

Based on extensive evidence, the combinatorialists Goulden and Jackson \cite[Conj.~1.2]{gjconj} had
conjectured that this combinatorial
problem had a surprising polynomial behavior: fixing $g$ and $n$, 
$H^g_{\al}$ is a simple combinatorial term times a symmetric
polynomial in $\al_1$, \dots, $\al_n$, with components in homogeneous
degree between $2g-3+n$ and $3g-3+n$.  This strongly suggests a
connection between this combinatorial problem and the moduli space of
curves!

Ekedahl, Lando, M. Shapiro, and Vainshtein explained this polynomiality
 with their ground-breaking ELSV-formula:
$$
H^g_{\al}
= r! \prod^n_{i=1}\left(  \frac {  \al_i^{\al_i}} {\al_i!} \right)
\int_{\cmbar_{g,n}} \frac { 1 - \la_1 + \cdots + (-1)^g \la_g}
{ (1-\al_1 \psi_1) \cdots (1- \al_n \psi_n)}
$$
Here the denominator should be considered formally inverted,
i.e.\ $1/(1-\al_i \psi_i) = 1 + \al_i \psi_i + \al_i^2 \psi_i^2 + \cdots$,
and the integral sign means to take the degree
of the  codimension $3g-3+n$ (dimension $0$) part of the integrand.
Expanding the integral in the ELSV formula yields
$$
\sum_{a_1 + \cdots + a_n +k=3g-3+n}
\left(  (-1)^k  \int_{\cmbar_{g,n}} \psi_1^{a_1} \cdots \psi_n^{a_n} \lambda_k \right)
( \al_1^{a_1} \cdots \al_n^{a_n} )
$$
and thus the polynomiality is explained by interpreting these numbers
as top intersections on the moduli space of curves!

Better yet, the highest-degree terms are precisely the subject of
Witten's conjecture, and indeed the five proofs of Witten's conjecture
using the ELSV formula involve in different ways getting at these
leading coefficients by asymptotic methods.  Furthermore,
the lowest-degree terms are precisely the subject of the $\la_g$-conjecture,
which provides the entree for the proof of \cite{lambdag}.

This polynomiality arises repeatedly, usually as a result of localization
on spaces of stable maps, and was a key ingredient in, for example,
work of Graber and the second author, e.g.\ \cite{gv,thmstar}.

By comparison, we give an example of polynomiality which has yet to be
satisfactorily explained geometrically, which allows us to clarify a 
vague but suggestive conjecture of \cite{gjv}. 

\subsection{A fourth sister?  Conjectural geometry behind polynomiality
of double Hurwitz numbers, \cite{gjv}, and work of Shadrin and Zvonkine}

We recall the history of the ELSV formula as motivation: we begin with
an enumerative problem in geometry, of Hurwitz numbers, which can also
be interpreted in terms of the combinatorics of the symmetric group.
For fixed $g$ and $n$, these numbers are (up to a combinatorial
factor) a polynomial in the parts of a partition (appearing in the
definition of Hurwitz number).  This polynomial is symmetric of degree
$3g-3+n$, which is highly suggestive of the moduli space of curves,
and indeed there is a genus $g$, $n$-pointed curve present in the
enumerative problem.  The generating function for these numbers is
constrained by the KP (Kadomtsev-Petviashvili) hierarchy; to see this,
a change of variables is necessary.  Thanks to the ELSV formula, these
numbers are related to the intersection theory of a fundamental moduli
space, and as a result, surprising structure is known on the
(tautological) cohomology ring of the moduli space.

The identical story will apply in the situation we now describe,
except that there is as of yet no moduli space, and no ELSV-type
formula.  There seems strong circumstantial evidence that there is a
moduli space completing the story; this is the content of Conjecture
3.5 of \cite{gjv}.

Instead of ``single'' Hurwitz numbers, we consider ``one-part double''
Hurwitz numbers $H^g_{(d),\be}$, defined in the same way as single Hurwitz numbers,
except that we require in addition complete branching over the point
$0 \in \proj^1$.    Then by means of character theory, 
\cite[Thm.~3.1]{gjv} shows that for fixed $g$ and $n$,  $H^g_{(d),\be}$
is a symmetric polynomial in $\be_1$, \dots, $\be_n$, whose
homogeneous pieces have even degree up to $4g-3+n$.

The fundamental question this suggests is: {\em is there a moduli
space and ELSV-type formula explaining this polynomiality?}

Note that there is a moduli space of dimension $4g-3+n$ ``present'' in
the problem: an $n$-pointed genus $g$ curve, along with a choice of
line bundle.  
Motivated by this, one possible answer is the following.

{\bf Conjecture \cite[Conj.~3.5]{gjv}.}
For each $g \geq 0$, $n \geq 1$, $(g,n) \neq (0,1), (0,2)$,
$$H^g_{(d),\be} = r^g_{(d), \be}! d \int_{\overline{\operatorname{Pic}}_{g,n}}
\frac { \Lambda_0 - \Lambda_2 + \cdots \pm \Lambda_{2g}}
{ (1- \be_1 \psi_1) \cdots (1 - \be_n \psi_n)}$$
where $\overline{\operatorname{Pic}}_{g,n}$, $\psi_i$, and $\Lambda_{2k}$
satisfy a number of reasonable properties.
For example,
$\overline{\operatorname{Pic}}_{g,n}$ should be a compactification
of the universal Picard variety over $\cm_{g,n}$, which supports a
(possibly virtual) fundamental class of the ``expected'' dimension $4g-3+n$.
$\psi_i$ should be an extension of the pullback of $\psi_i$ from $\cm_{g,n}$.
See \cite{gjv} for a complete list of proposed properties.

It should be emphasized that the fundamental question is that of
finding the right space, and this conjecture should be seen as merely
a proposal.  Some of the properties suggested by geometry yield
testable constraints on double Hurwitz numbers (such as the string and
dilaton equation), and these indeed hold \cite[Prop.~3.10]{gjv}.

Shadrin and Zvonkine show much more in \cite{sz} and \cite{s}. We
emphasize that independent of the conjecture, their results are
meaningful statements about double Hurwitz numbers. But their
unpacking of the structure of these Hurwitz numbers is predicted
precisely by the form of the conjecture. We state their results
imprecisely in order to emphasize their form. Let $H$ be the
generating function for the double Hurwitz numbers
$H^g_{(d),\be}$. Let $U$ be the generating function for just the
highest-degree terms in the conjecture (i.e.\ those with no
$\Lambda_i$ terms for $i>0$). Shadrin and Zvonkine define a linear
differential operator $L$ which arises naturally in the cut-and-join
formalism ($L = \sum b p_b \frac {\partial} {\partial p_b}$, where
$p_b$ tracks parts of $\al$ of size $b$).

The KP hierarchy is a system of partial differential equations on a power series $F$
in an infinite set of variables $p_1$, $p_2$, \dots.  The first three equations are:
\begin{eqnarray*}
F_{2,2} &=& - \frac 1 2 F^2_{1,1} + F_{3,1} - \frac 1 {12} F_{1,1,1,1} \\
F_{3,2} &=& - F_{1,1} F_{2,1} + F_{4,1} - \frac 1 6 F_{2,1,1,1} \\
F_{4,2} &=& - \frac 1 2 F^2_{2,1} - F_{1,1} F_{3,1} + F_{5,1} + \frac 1 8 F^2_{1,1,1}
+ \frac 1 {12} F_{1,1} F_{1,1,1,1} \\
& & \quad \quad  - \frac 1 4 F_{3,1,1,1} + \frac 1 {120} F_{1,1,1,1,1,1}
\end{eqnarray*}
Subscript $i$ refers to differentiation by $p_i$.  The exponent $\tau = e^F$ 
of any solution is called a $\tau$-function of the hierarchy.
For readable expositions on the KP hierarchy, see \cite[\S 6]{kazarian}
and (influenced by this) \cite{gjkp}.

Shadrin and Zvonkine prove that under a scaling and renaming of
variables, {\bf (i)} $L^2 H$ is a $\tau$-function of the KP hierarchy,
i.e.\ it satisfies the bilinear Hirota equations.  Furthermore, $L^2
H$ satisfies the linearized KP equations.  {\bf (ii)} $U$ is a
$\tau$-function for the KP hierarchy (in unusual variables), and thus
satisfies the bilinear Hirota equations.  Furthermore, it satisfies
the linearized KP equations in the same variables.  This is a complete
analogue of Witten's Conjecture.  Statement {\bf (i)} follows in a
standard way from the general theory of integrable systems, but the
proof of {\bf (ii)} is quite subtle.  

Shadrin has recently taken this further.  In \cite{s}, he applies Kazarian's
techniques to this series, relating the conjectured ``intersection
numbers'' to the Hirota equations.  The computations turn out to be
simpler, and he gives explicit and rather transparent formulas for the
generating series.

One might take Shadrin and Zvonkine's results to suggest that the
moduli space behind the integrable system is precisely that suggested
by \cite[Conj.~3.5]{gjv}. However, this may not be the case: their
methods are quite robust, and similar spaces (for example, with $n+1$
marked points, with the last point the preimage of $0$) should yield
similar results.

\section{Teleman's work on family topological field theory} \label{s:teleman}

\subsection{Terminology} 
For this section only, dimensions here mean dimensions over $\mathbb{R}$ and
surfaces mean topological oriented 2-$\dim_{\mathbb{R}}$ 
``manifolds with nodal singularities''.
Curves, however, mean algebraic (or complex) curves!

\subsection{Topological field theory and moduli of curves}
C.~Teleman introduced the notion of 
\emph{family topological field theory} (FTFT) into Gromov--Witten theory.
In the semisimple case, he completely classified all FTFTs.
As an important application, he proved \emph{Givental's conjecture} \cite{aG1}
on the quantization formula, which in particular gives an explicit
reconstruction of higher genus semisimple Gromov--Witten theory from the genus
zero theory.

Teleman's result is 
the first instance of applying two powerful results from \emph{topology} 
of moduli spaces of curves to Gromov--Witten theory: \emph{Harer stability}
and the \emph{Madsen--Weiss theorem} (Mumford's conjecture).
Harer stability is purely topological in its formulation, 
and Mumford's conjecture so far has no algebro-geometric proof.
In this section, we will explain how these two results can be applied to 
Gromov--Witten theory.

\subsection{Two dimensional topological field theory and semisimplicity}
A 2-$\dim_{\mathbb{R}}$ topological field theory (TFT) is a symmetric,
(strong, monoidal) functor of 
topological tensor categories $Z: \mathbf{C} \to \mathbf{Vec}$.
$\mathbf{Vec}$ is the usual tensor category of vector spaces.
$\mathbf{C}$ is the category whose objects are 1-$\dim_{\mathbb{R}}$ 
oriented closed manifolds, i.e.~a disjoint union of oriented $S^1$'s.
The morphisms are oriented cobordisms of the objects, with obvious 
compositions. The tensor structure is defined by disjoint union of objects.

The notion of TFT is equivalent to the notion of Frobenius algebra. 
A (commutative) \emph{Frobenius algebra} is a $k$-algebra with an identity $1$,
a pairing $(\cdot, \cdot) : A^{\otimes 2} \to A$ which is symmetric and 
nondegenerate, and satisfies the \emph{Frobenius property}
\[
 (a * b, c) = (a, b*c),
\]
where $*$ stands for the multiplicative structure of $A$.

The equivalence of TFT and Frobenius algebra can be seen as follows. 
$A = Z (S^1)$, and the identity element $1 \in A$ is defined by
$Z(\text{cap})$,  (0 input, 1 output).
The nondegenerate pairing $(\cdot, \cdot)$ is defined by 
$Z(\text{bent cylinder})$ (2 inputs, 0 output), and the ring structure is 
defined by $Z(\text{pair of pants})$ (2 inputs, 1 output).

\vspace{10pt}

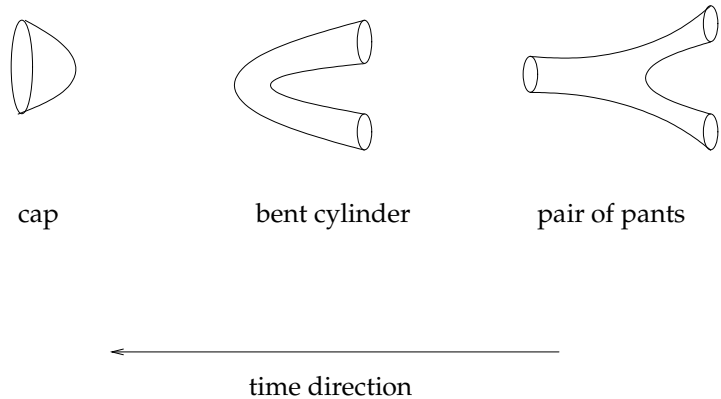
\begin{center}

\setlength{\unitlength}{0.0005in}
\begingroup\makeatletter\ifx\SetFigFont\undefined%
\gdef\SetFigFont#1#2#3#4#5{%
  \reset@font\fontsize{#1}{#2pt}%
  \fontfamily{#3}\fontseries{#4}\fontshape{#5}%
  \selectfont}%
\fi\endgroup%
{\renewcommand{\dashlinestretch}{30}
\begin{picture}(7366,4092)(0,-10)
\put(120,3428){\ellipse{224}{974}}
\put(3683,3690){\ellipse{150}{450}}
\put(3683,2753){\ellipse{150}{374}}
\put(7283,3878){\ellipse{150}{374}}
\put(7283,2752){\ellipse{150}{374}}
\put(5408,3353){\ellipse{150}{374}}
\path(5708,465)(1058,465)
\path(1178.000,495.000)(1058.000,465.000)(1178.000,435.000)
\path(158,3915)(160,3914)(164,3912)
	(171,3908)(182,3902)(197,3894)
	(217,3883)(240,3870)(267,3855)
	(296,3838)(328,3819)(361,3798)
	(395,3777)(429,3754)(462,3730)
	(495,3706)(526,3680)(555,3653)
	(583,3626)(608,3597)(630,3566)
	(650,3534)(665,3500)(676,3464)
	(683,3428)(683,3390)(677,3353)
	(665,3317)(648,3284)(628,3252)
	(604,3223)(577,3196)(548,3170)
	(516,3146)(483,3124)(448,3103)
	(412,3083)(376,3064)(338,3045)
	(302,3028)(266,3013)(231,2998)
	(200,2985)(171,2973)(146,2964)
	(125,2956)(109,2949)(97,2945)
	(89,2942)(85,2941)(83,2940)
\path(3683,3915)(3682,3915)(3678,3914)
	(3672,3912)(3662,3910)(3648,3906)
	(3630,3901)(3607,3895)(3579,3888)
	(3546,3879)(3508,3869)(3466,3857)
	(3420,3845)(3371,3831)(3319,3816)
	(3264,3799)(3208,3783)(3150,3765)
	(3092,3746)(3034,3727)(2976,3707)
	(2919,3687)(2862,3666)(2807,3644)
	(2754,3621)(2702,3598)(2653,3574)
	(2605,3550)(2561,3524)(2519,3497)
	(2480,3469)(2445,3440)(2413,3409)
	(2386,3377)(2364,3344)(2347,3310)
	(2337,3275)(2333,3240)(2337,3205)
	(2347,3170)(2364,3136)(2386,3103)
	(2413,3071)(2445,3040)(2480,3011)
	(2519,2983)(2561,2956)(2605,2930)
	(2653,2906)(2702,2882)(2754,2859)
	(2807,2836)(2862,2814)(2919,2793)
	(2976,2773)(3034,2753)(3092,2734)
	(3150,2715)(3208,2697)(3264,2681)
	(3319,2664)(3371,2649)(3420,2635)
	(3466,2623)(3508,2611)(3546,2601)
	(3579,2592)(3607,2585)(3630,2579)
	(3648,2574)(3662,2570)(3672,2568)
	(3678,2566)(3682,2565)(3683,2565)
\path(3683,3465)(3681,3465)(3676,3464)
	(3667,3463)(3654,3462)(3635,3460)
	(3611,3458)(3580,3454)(3545,3451)
	(3505,3446)(3460,3441)(3412,3436)
	(3362,3429)(3310,3423)(3256,3416)
	(3203,3408)(3150,3400)(3098,3392)
	(3048,3383)(2999,3374)(2953,3364)
	(2909,3354)(2868,3343)(2831,3331)
	(2797,3318)(2767,3304)(2743,3290)
	(2724,3274)(2712,3257)(2708,3240)
	(2712,3222)(2724,3205)(2743,3188)
	(2767,3171)(2797,3155)(2831,3140)
	(2868,3125)(2909,3111)(2953,3097)
	(2999,3083)(3048,3070)(3098,3058)
	(3150,3046)(3203,3034)(3256,3022)
	(3310,3011)(3362,3000)(3412,2990)
	(3460,2981)(3505,2973)(3545,2965)
	(3580,2958)(3611,2953)(3635,2949)
	(3654,2945)(3667,2943)(3676,2941)
	(3681,2940)(3683,2940)
\path(5408,3540)(5409,3540)(5411,3540)
	(5416,3539)(5423,3539)(5433,3538)
	(5446,3537)(5462,3536)(5481,3534)
	(5504,3532)(5529,3531)(5557,3529)
	(5588,3527)(5621,3526)(5655,3525)
	(5692,3524)(5729,3523)(5768,3523)
	(5808,3523)(5849,3524)(5891,3525)
	(5934,3527)(5979,3530)(6025,3534)
	(6073,3539)(6123,3546)(6174,3553)
	(6228,3562)(6283,3573)(6341,3585)
	(6399,3599)(6458,3615)(6520,3633)
	(6579,3653)(6636,3673)(6689,3693)
	(6739,3714)(6785,3734)(6829,3755)
	(6869,3775)(6908,3796)(6944,3816)
	(6978,3837)(7010,3857)(7041,3878)
	(7070,3898)(7098,3917)(7124,3937)
	(7149,3955)(7172,3973)(7193,3990)
	(7212,4005)(7229,4019)(7243,4031)
	(7255,4041)(7265,4049)(7272,4055)
	(7277,4060)(7280,4063)(7282,4064)(7283,4065)
\path(5408,3165)(5409,3165)(5411,3165)
	(5416,3165)(5423,3165)(5433,3165)
	(5446,3165)(5462,3164)(5481,3164)
	(5504,3164)(5529,3163)(5557,3162)
	(5588,3161)(5621,3159)(5655,3158)
	(5692,3156)(5729,3153)(5768,3150)
	(5808,3147)(5849,3143)(5891,3138)
	(5934,3133)(5979,3126)(6025,3119)
	(6073,3111)(6123,3101)(6174,3091)
	(6228,3079)(6283,3065)(6341,3050)
	(6399,3033)(6458,3015)(6520,2994)
	(6579,2973)(6636,2951)(6689,2930)
	(6739,2908)(6785,2887)(6829,2866)
	(6869,2845)(6908,2824)(6944,2804)
	(6978,2784)(7010,2764)(7041,2744)
	(7070,2725)(7098,2706)(7124,2687)
	(7149,2669)(7172,2652)(7193,2636)
	(7212,2622)(7229,2609)(7243,2597)
	(7255,2588)(7265,2580)(7272,2574)
	(7277,2570)(7280,2567)(7282,2566)(7283,2565)
\path(7283,3690)(7281,3689)(7276,3688)
	(7267,3685)(7253,3681)(7234,3676)
	(7210,3669)(7181,3660)(7148,3650)
	(7111,3638)(7071,3625)(7030,3611)
	(6988,3596)(6945,3581)(6903,3564)
	(6862,3547)(6822,3529)(6785,3510)
	(6750,3490)(6717,3469)(6687,3447)
	(6661,3423)(6639,3398)(6623,3372)
	(6612,3344)(6608,3315)(6612,3286)
	(6623,3258)(6639,3232)(6661,3207)
	(6687,3183)(6717,3161)(6750,3140)
	(6785,3120)(6822,3101)(6862,3083)
	(6903,3066)(6945,3049)(6988,3034)
	(7030,3019)(7071,3005)(7111,2992)
	(7148,2980)(7181,2970)(7210,2961)
	(7234,2954)(7253,2949)(7267,2945)
	(7276,2942)(7281,2941)(7283,2940)
\put(83,1815){\makebox(0,0)[lb]{\smash{{\SetFigFont{10}{14.4}{\rmdefault}{\mddefault}{\updefault}cap}}}}
\put(2558,1815){\makebox(0,0)[lb]{\smash{{\SetFigFont{10}{14.4}{\rmdefault}{\mddefault}{\updefault}bent cylinder}}}}
\put(5483,1815){\makebox(0,0)[lb]{\smash{{\SetFigFont{10}{14.4}{\rmdefault}{\mddefault}{\updefault}pair of pants}}}}
\put(2483,15){\makebox(0,0)[lb]{\smash{{\SetFigFont{10}{14.4}{\rmdefault}{\mddefault}{\updefault}time direction}}}}
\end{picture}
}

\end{center}

\vspace{10pt}

We will mainly be concerned with the \emph{semisimple} case, where
the product structure is diagonalizable.
That is, $A \cong \oplus_i k \epsilon_i$ as an $k$-algebra, 
where $\epsilon_i * \epsilon_j = \delta_{ij} \epsilon_i$ is the 
\emph{canonical basis}.
Therefore, up to isomorphisms, $A$ is classified by 
$\Delta_i := (\epsilon_i, \epsilon_i)$.
Let 
\[
 \tilde{\epsilon}_i := \frac{1}{\sqrt{\Delta_i}} \epsilon_i
\]
be the \emph{normalized canonical basis}.

\begin{proposition} \label{p:1}
Let ${}^qC_g^p$ be an oriented surface of genus $g$, with $p$ inputs and
$q$ outputs.
Then $Z({}^qC_g^p)$ defines a map 
\[
 {}^q\mu_g^p: A^{\otimes p} \to A^{\otimes q},
\]
such that ${}^q\mu_g^p$ is ``diagonal'' in the (tensor power of the)
normalized canonical basis.
Furthermore, the entry of ${}^q\mu_g^p$ in 
$\tilde{\epsilon}_i^{\otimes p} \mapsto \tilde{\epsilon}_i^{\otimes q}$ is
$\Delta_i^{\chi({}^qC_g^p)/2}$, where $\chi({}^qC_g^p)$ is the (topological) 
Euler characteristic of ${}^qC_g^p$.
\end{proposition}

This proposition can be easily proved by decomposing ${}^qC_g^p$ into
pairs of pants and applying the facts that 
$Z({}^1C_0^1)$ is the identity and that $Z({}^1C_0^2)$ defines the algebra
multiplication and is diagonal with entries $\Delta_i^{1/2}$ in the
normalized canonical basis. 

Let ${}^1C_1^1$ be the torus with one outgoing and one incoming boundary.
Since the ring structure is determined by $\tilde{Z}$ on a pair of pants,
one can piece together 2 pairs of pants and obtain ${}^1C_1^1$.
For future references, denote $\alpha$ (the diagonal matrix) 
${}^1\mu_1^1=\tilde{Z}({}^1C_1^1)$.

\subsection{Family topological field theories}
A FTFT assigns a family of topological $(m,n)$-pointed
2-$\dim_{\mathbb{R}}$ surfaces $\mathcal{C} \to B$ a cohomology class 
$\bar{Z}(B) \in H^*(B, \operatorname{Hom}(A^m, A^n))$,
where $A^m$ are, in general, local systems on $B$.
$\bar{Z}$ must satisfy two additional properties, in addition to being a
``fiberwise TFT:'' 
\begin{enumerate}
 \item functoriality with respect to the (topological) base change;
 \item a strong gluing axiom.
\end{enumerate}
The first property ensures that all FTFTs are pullbacks of
the universal setting over the classifying spaces.
The second property requires a little more explanation.

The gluing requires parameterizations of boundary circles, instead of
the punctures.
The strong gluing axiom asserts that the gluing axiom has to hold for 
\emph{any} lifting to families with parameterized boundaries.
The $(m,n)$-pointed family can be lifted to a family of $m$ incoming and 
$n$ outgoing circles in the following way.
Take the torus bundle $\tilde{B} \to B$ with fiber $(S^1)^m \times (S^1)^n$
the product of unit tangent spaces at the marked points.
(Since the space of Riemannian metrics is contractible, the choice
of the metric does not matter.)
Because $\operatorname{Diff}_+(S^1)$, the orientation-preserving 
diffeomorphisms, is homotopy equivalent to its subgroup
of rigid rotations, $\tilde{B}$ is homotopy equivalent to the base space of 
the corresponding family with parameterized circle boundaries.
Note that $\mathcal{C}$ allows nodal singularities.

Teleman classified all FTFTs for which the fiberwise TFT is
generically semisimple.

Teleman established this classification in 3 steps.

\textbf{Step 1.} 
Classification of FTFTs $\tilde{Z}$ associated to smooth families of surfaces 
with parameterized boundaries. 

\textbf{Step 2.}
Going from parameterized boundaries to punctures. This requires a new piece
of information, which is denoted $E(\psi)$.
Roughly, at each outgoing boundary circle, $E(\psi)$ 
(or $E(\psi)^{-1}$ for an incoming circle) ``transforms'' the parameterized 
boundary to a puncture.

\textbf{Step 3.}
Allowing nodal (stable) degenerations of curves. 
The corresponding FTFT is denoted $\bar{Z}$.
$\bar{Z}$ requires again one new piece of data: 
$L:=\bar{Z}(\text{pinched cylinder})$. 
In fact, this classification can be generalized to, what Teleman calls
``Lefschetz theory'', where the family of curves does not have to be stable,
and with more general gluing axioms.
Interested readers are referred to Teleman's original article \cite{cT}.

\vspace{10pt}

\begin{center}

\setlength{\unitlength}{0.0006in}
\begingroup\makeatletter\ifx\SetFigFont\undefined%
\gdef\SetFigFont#1#2#3#4#5{%
  \reset@font\fontsize{#1}{#2pt}%
  \fontfamily{#3}\fontseries{#4}\fontshape{#5}%
  \selectfont}%
\fi\endgroup%
{\renewcommand{\dashlinestretch}{30}
\begin{picture}(2716,1971)(0,-10)
\put(158,1382){\ellipse{300}{1124}}
\put(2558,1382){\ellipse{300}{1124}}
\path(158,1944)(160,1944)(163,1943)
	(169,1941)(180,1938)(194,1934)
	(213,1929)(236,1922)(265,1914)
	(298,1905)(336,1894)(378,1882)
	(423,1869)(471,1854)(522,1839)
	(575,1822)(628,1805)(683,1787)
	(738,1768)(792,1749)(846,1729)
	(898,1709)(950,1688)(999,1667)
	(1047,1645)(1093,1622)(1136,1598)
	(1177,1574)(1214,1548)(1249,1522)
	(1279,1495)(1306,1466)(1328,1437)
	(1344,1406)(1354,1375)(1358,1344)
	(1354,1313)(1344,1283)(1328,1253)
	(1306,1225)(1279,1198)(1249,1173)
	(1214,1149)(1177,1126)(1136,1104)
	(1093,1083)(1047,1063)(999,1044)
	(950,1025)(898,1008)(846,991)
	(792,974)(738,958)(683,943)
	(628,929)(575,915)(522,902)
	(471,889)(423,878)(378,867)
	(336,858)(298,849)(265,842)
	(236,836)(213,831)(194,826)
	(180,823)(169,821)(163,820)
	(160,819)(158,819)
\path(2558,1944)(2556,1944)(2553,1943)
	(2547,1941)(2536,1938)(2522,1934)
	(2503,1929)(2480,1922)(2451,1914)
	(2418,1905)(2380,1894)(2338,1882)
	(2293,1869)(2245,1854)(2194,1839)
	(2141,1822)(2088,1805)(2033,1787)
	(1978,1768)(1924,1749)(1870,1729)
	(1818,1709)(1766,1688)(1717,1667)
	(1669,1645)(1623,1622)(1580,1598)
	(1539,1574)(1502,1548)(1467,1522)
	(1437,1495)(1410,1466)(1388,1437)
	(1372,1406)(1362,1375)(1358,1344)
	(1362,1313)(1372,1283)(1388,1253)
	(1410,1225)(1437,1198)(1467,1173)
	(1502,1149)(1539,1126)(1580,1104)
	(1623,1083)(1669,1063)(1717,1044)
	(1766,1025)(1818,1008)(1870,991)
	(1924,974)(1978,958)(2033,943)
	(2088,929)(2141,915)(2194,902)
	(2245,889)(2293,878)(2338,867)
	(2380,858)(2418,849)(2451,842)
	(2480,836)(2503,831)(2522,826)
	(2536,823)(2547,821)(2553,820)
	(2556,819)(2558,819)
\put(383,69){\makebox(0,0)[lb]{\smash{{\SetFigFont{10}{14.4}{\rmdefault}{\mddefault}{\updefault}pinched cylinder}}}}
\end{picture}
}

\end{center}

\vspace{10pt}

Teleman proved that $\tilde{Z}$, $E(\psi)$, and $L$
are sufficient to ``piece together'' $\bar{Z}$.
Roughly, one can replace punctures by parameterized boundaries via $E(\psi)$,
and smooth the nodes via $L$.
The remaining theory is then $\tilde{Z}$.
Among the 3 steps, Step 1 involves the aforementioned Harer stability and 
Madsen--Weiss Theorem and is most relevant to this article's audience.
In the remaining of this section, we will give some explanation of 
Step 1.
But we first give a short digression on how Teleman's classification
proves Givental's conjecture.

\subsection{Proof of Givental's conjecture}
We will use the following terminology proposed in \cite{ypL4}.

A \emph{geometric} Gromov--Witten theory (GWT) is the ``usual'' GWT
constructed from moduli spaces of stable maps of \emph{all genera}.
Each geometric GWT is encoded in a single generating function
\begin{equation} \label{e:tauGW}
 \tau_{GW} := \exp \left( \sum_{g=0}^{\infty} \hbar^{g-1} F_g \right)
\end{equation}
where $F_g$ is the genus $g$ generating function of GW invariants.

An \emph{axiomatic} GWT is the one constructed from Givental's
quantization formalism, when the Frobenius structure is semisimple.
Roughly, the generating function of an axiomatic GWT is \emph{defined}
to be
\begin{equation} \label{e:tauG}
  \tau_G := \hat{\operatorname{O}} \tau^{N pt},
\end{equation}
where $\tau^{N pt}$ is the $N$ copies of Witten--Kontsevich generating 
function (or equivalently GWT of $N$ points), and $\hat{\operatorname{O}}$
is an operator defined by quantizing a quadratic functions defined via
\emph{genus zero} GW invariants (or equivalently the Frobenius structure).

When the Frobenius structure is semisimple, one can ask whether
$\tau_G$ is equal to $\tau_{GW}$.
The equality is termed \emph{Givental's conjecture}.
The conjecture in particular implies that all higher genus GW invariants
can be reconstructed from genus zero invariants, as $\hat{O}$ involves
only genus zero data.

One can prove that $\tau_G$ can be written as
\[
 \tau_G = \exp \left( \sum_{g=0}^{\infty} \hbar^{g-1} G_g \right).
\]
Furthermore, it can be shown without much difficulty that $G_0=F_0$, 
which implies that the genus zero truncation of Givental's conjecture holds.
See \cite[Part II]{LP}  for details.

Teleman's classification implies Givental's conjecture.
Unfortunately, the brief account given here is mostly redundant for those 
who are familiar with the basic definitions.
A meaningful explanation of this involves the detailed construction of
Givental's quantization formalism and is beyond the scope of this article.
Interested readers are referred to Teleman's original paper \cite{cT}.

First of all, it is easy to see that axiomatic GWTs are FTFT.
For those who are familiar with the definition of \emph{cohomological field
theory}, GWTs are CohFTs and CohFTs satisfy the axioms of FTFTs.
In fact, Teleman's classification implies a \emph{stronger} versions of 
Givental's conjecture: 
The corresponding axiomatic and geometric CohFTs are equal.
Indeed, CohFTs yield cohomology classes on $\Mbar_{g,n}$, and the
integration of these classes (with monomials of $\psi$-classes) gives
GW invariants.

To prove Givental's conjecture, Teleman traces through Givental's quantization
formalism and shows that all semisimple GWTs can be reconstructed from
the genus zero data, and from known facts about the tautological classes
of the moduli space of curves.
(Very naively speaking, one can see all discussion in this section 
involves only genus zero data and moduli of curves.)
Since axiomatic and geometric theories are, by construction, 
identical in genus zero, the proof is then complete.

\subsection{Wheeled PROP structure} 
Let ${}^q{M}_g^p$ be the moduli space of the Riemann surfaces
of genus $g$ with $p$ incoming and $q$ outgoing marked points.
Let ${}^q\tilde{M}_g^p$ be the torus bundle over ${}^q{M}_g^p$ as discussed 
earlier.
Homotopically, ${}^q{M}_g^p \sim B \Gamma_g^{p+q}$, where $\Gamma_g^{p+q}$ is
the corresponding mapping class group, and one can think of ${}^q{M}_g^p$ as 
a classifying space of bordered Riemann surfaces of genus $g$ with $p$ 
incoming and $q$ outgoing parameterized boundaries.
Gluing the outgoing boundary circles of a surface with incoming ones of another
gives the structure of a \emph{PROP} on ${}^q\tilde{M}_g^p$.
Roughly, a PROP is something like an operad, but with multiple outputs.\footnote{The formal definition of PROP, or product and permutation category,
can be found in, e.g.~\cite{mM}. Briefly, a PROP
$(P,*,S, 1))$ is a symmetric strict monoidal category such that
\begin{enumerate}
\item the objects $Ob(P)$ are identified with the set $\mathbb{Z}_{\ge 0}$;
\item the product satisfies $m*n = m+n$, for any $m,n \in Ob(P)$ 
(hence the unit $1=0$);
\item $S$ is the permutation symmetry;
\item each hom-set $Mor_P(m,n)$ is a $k$-module and the operations of the
monoidal category $P$ are compatible with this $k$-linear structure.
\end{enumerate}}
In addition to the PROP structure, it also allows \emph{contraction}, via
gluing boundary circles of the same surface. 
The additional contraction operation makes it a \emph{wheeled PROP} \cite{MMS}.

\subsection{Harer stability}
We recall Harer stability. 
For the notational convenience, we will \emph{not} distinguish the outgoing 
and incoming boundary, but use $r := p+q$.
Let $\tilde{M}_{g,r}$ be the (homotopy type of the) moduli space of 
oriented surfaces of genus $g$, and $r$ boundary circles.
(The number of punctures will be fixed during the discussion, and will be 
omitted from the notations.)
Let $C_{g,r}$ be a bordered oriented surface of genus $g$ with $r$ boundary
circles. 
Define three operations
 \begin{align*}
  \phi_1 &: C_{g,r} \to C_{g+1, r-1},  &r \ge 2, \\
  \phi_2 &: C_{g,r} \to C_{g, r+1},    &r \ge 1, \\
  \phi_3 &: C_{g,r} \to C_{g+1, r-2},  &r \ge 2
 \end{align*}
as follows:
$\phi_1$ is defined by gluing a pair of pants along 2 boundary circles;
$\phi_2$ is defined by gluing a pair of pants along 1 boundary circle;
$\phi_3$ is defined by gluing 2 boundary circles together.

\begin{theorem}[Harer stability] \label{t:harer}
$(\phi_1)_* : H_k (\tilde{M}_g^r) \to H_k (\tilde{M}_{g+1}^{r-1})$,
and similarly $(\phi_2)_*$ and $(\phi_3)_*$, are isomorphisms for
$k \le [g/3]$ \emph{(the stable range)}.
\end{theorem}

The stable range can be improved, but we only need this estimate.

Thanks to Harer stability, one can talk about $H_k({M}_{\infty})$.
Here $r$ is redundant in the stable range by a combination of $\phi_i$'s,
and hence $\tilde{M}_{\infty} = M_{\infty}$.

\subsection{Classification of ${}^1\tilde{Z}$}
We will first study ${}^1\tilde{Z}$, i.e.~FTFT with 
one outgoing and no incoming boundary.
One can increase genus by applying the composition $\phi_2 \circ \phi_1$.
By semisimplicity, the effect of this operation on ${}^1\tilde{Z}_g$ 
is the multiplication by (the diagonal matrix) $\alpha$. 
Therefore, ${}^1\tilde{Z}_{g+G} = (\alpha^G) {}^1\tilde{Z}_g$ 
in the stable range and $(\alpha^{-G}) {}^1\tilde{Z}_G$ 
stabilizes to $\tilde{Z}^+ \in H^*(M_{\infty})$.

\begin{proposition} \label{p:2}
(i) The degree zero component of $\tilde{Z}^+ \in H^*(M_{\infty})$ is the 
identity.

(ii) $H^*(M_{\infty})$ carries a Hopf algebra structure.

(iii) $\log \tilde{Z}^+ =\sum_{l \ge 1} a_l \kappa_l$, where $a_l \in A$.
\end{proposition}

(i) can be seen by restricting our family of oriented smooth surfaces to
a point, and apply base change property:
For an inclusion $\iota: pt \to M$, the pullback in cohomology in degree zero
$\iota^* :H^0(M) \to H^0(pt)$ induces isomorphism on degree zero
component of $\tilde{Z}^+$.

(ii) is a consequence of the wheeled PROP structure.

(iii) is a consequence of the wheeled PROP structure, gluing axiom, and
Mumford's conjecture. 
We will explain the idea now.

An element in a Hopf algebra is called {\em primitive} if 
$\Delta (x) = 1 \otimes x + x \otimes 1$.
It is {\em group-like} if $\Delta (g) = g \otimes g$.
It is easy to check that group-like elements are exponentials of 
primitive elements.

\begin{lemma}
$\tilde{Z}^+$ is a group-like element.
\end{lemma}

This follows from the monoidal structure (in PROP) on $M_{\infty}$
defined by gluing two surfaces into a pair of pants:
\[
 \tilde{M}^1_{g_1} \times \tilde{M}^1_{g_2} \to \tilde{M}^1_{g_1+ g_2}.
\]
Note that the multiplicative factor $\alpha^{-(g_1+g_2)}$ is consistent
on both sides.
Therefore, $\log \tilde{Z}^+$ is a primitive element.
Now, Madsen--Weiss's theorem:

\begin{theorem}[Mumford's conjecture, \cite{MW}] 
\[
 H^*(M_{\infty}, \mathbb{Q}) = \mathbb{Q} [\kappa_0, \kappa_1, \kappa_2, ...].
\]
\end{theorem}

Thanks to the polynomial structure of the Hopf algebra, the only primitive
elements are $\sum_{l \ge 0} a_l \kappa_l$.
Applying Proposition~\ref{p:2} (i), $a_0=0$ and
$\log \tilde{Z}^+ =\sum_{l \ge 1} a_l \kappa_l$.

Note that the converse of the statement also holds: All elements of the form
$\sum_{l \ge 1} a_l \kappa_l$ serve as FTFT on smooth surfaces with 
parameterized boundaries.
This completely classifies $\log \tilde{Z}^+$ and hence ${}^1\tilde{Z}$.

\subsection{The final steps:  Classifying ${}^q\tilde{Z}^p$}
One can increase genus by gluing ${}^1C^1$, and can increase $p$ and $q$ by
gluing a genus zero surface with $p+1$ inputs and $q$ outputs to the existing
one output.
The latter gives a map from ${}^1\tilde{M}_g$ to ${}^q\tilde{M}^p_g$.
Thanks to Harer stability, both operations are homology equivalences 
in the stable range.
By  functoriality of base change, ${}^q\tilde{Z}_g^p$ is determined by
its pullback to ${}^1\tilde{M}_g$.
On ${}^1\tilde{M}_g$, this can be seen as feeding the output of 
${}^1\tilde{M}_g$ into ${}^q\mu^{p+1}: A^{\otimes p+1} \to A^{\otimes q}$.
Therefore:

\begin{corollary}
${}^q \tilde{Z}_g^p$ is diagonal in the (tensor power of the) normalized
basis, and the entry for 
$\tilde{\epsilon}_i^{\otimes p} \mapsto \tilde{\epsilon}_i^{\otimes q}$ is
$\Delta_i^{\chi/2} \exp \left( \sum_{l \ge 1} a_{il} \kappa_l \right)$.
\end{corollary}
This completes the classification of $\tilde{Z}$.

\section{Witten's conjecture on $r$-spin curves} \label{s:invariance}
In this section, we  briefly describe 
Faber--Shadrin--Zvonkine's proof \cite{FSZ} 
 of  Witten's conjecture on $r$-spin curves.
This conjecture also follows from combining Teleman's more general result 
in the previous section, and Givental's result in \cite{aG2}.
\cite{FSZ} takes a different approach.

Semisimplicity of the Frobenius structure 
is tacitly assumed throughout this section.

\subsection{Overview}
Witten's conjecture states that a certain kind of Gromov--Witten theory
can be constructed from $r$-spin curves, whose moduli spaces are
certain branched covers of $\Mbar_{g,n}$.
Furthermore, the generating function $\tau_W$ 
constructed from Witten's correlators satisfies the Gelfand--Dickii hierarchy.

The proof has several ingredients.
Here we explain \emph{one} aspect of  the proof,
which is closer to our own.

\textbf{Step 0.}
Witten \cite{witten2} and Jarvis--Kimura--Vaintrob \cite{JKV} show,
building on Jarvis' earlier work \cite{tJ1, tJ2},
that the genus zero truncation of the conjecture holds.
Several rigorous constructions of Witten's correlators for all genera 
were defined (e.g.~\cite{tM, PV}).

\textbf{Step 1.}
The first author proposed to prove Witten's conjecture via a version of
Givental's conjecture:
First, show that the axiomatic generating function $\tau_G$ satisfies
Gelfand--Dickii hierarchy.
Then, show 
\begin{equation} \label{e:WG}
 \tau_W = \tau_G, 
\end{equation}
where $\tau_W$ is Witten's generating function.

\textbf{Step 2.}
Givental showed that $\tau_G$ satisfies Gelfand--Dickii hierarchy \cite{aG2}.

\textbf{Step 3.}
The first author proved a $g \le 2$ truncation of \eqref{e:WG} 
by introducing an algorithm, termed \emph{invariance constraints}, 
to compute universal relations in axiomatic GWT.
The \emph{universal relations} are, be definition,
the relations which hold for all axiomatic GWTs.
Furthermore, it is shown that the tautological relations are universal
relations \cite{ypL1, ypL2, ypL3}.

\textbf{Step 4.}
Faber--Shadrin--Zvonkine proved \eqref{e:WG} in full generality by
showing a reconstruction theorem which allows one to reduce the comparison
of the geometric GWT and the axiomatic GWT to a genus zero version.
It can be considered as a special case of Teleman's theorem
in the previous section.
This proof uses the idea of invariance constraints and Ionel's theorem
\cite{eI}, and is entirely different from Teleman's approach.

It is worth noting that Witten's conjecture is related to the Frobenius
structure of type $A_r$ singularities.
A generalization to type $D$ and $E$ has been carried out by H.~Fan,
T.~Jarvis, and Y.~Ruan (in preparation).

We will briefly explain Steps 0, 3 and 4.

\subsection{Witten's correlators}
\label{wc}An $r$-spin structure on a smooth curve $C$ of genus $g$ 
with $n$ marked points $x_i$
is a line bundle $L$ together with the identification
\[
 L^{\otimes r} = K_C \left( \sum_i m_i x_i \right),
\]
where $0 \le m_i \le r-1$ are integers such that $r | (2g-2-\sum_i m_i)$.
The moduli space of $r$-spin structures has a natural compactification
$\Mbar_{g,m}^{1/r}$ with a natural morphism $\pi$ to 
$\Mbar_{g,n}$, which is a degree $r^{2g-1}$ ``branched cover''.
Note that there is a global automorphism (multiplication by $r$th roots of 
unity) of any point in $\Mbar^{1/r}_{g,n}$.  
The degree is $r^{2g}$ for the corresponding coarse moduli spaces.

Witten's correlators are defined by a virtual fundamental class construction
\cite{tM, PV, aC} on $\Mbar_{g,(m_1,...,m_n)}^{1/r}$
\[
 \int_{[\Mbar_{g,(m_1,...,m_n)}^{1/r}]^{\vir}} \pi^* \Psi,
\]
where $\Psi$ is a monomial of $\psi$-classes on $\Mbar_{g,n}$.
By the projection formula, it is equal to 
$\int_{[\Mbar_{g,n}]} c_W \Psi$, where
$c_W (m) \cap \Mbar_{g,n}:= \pi_*([\Mbar_{g,a}^{1/r}]^{\vir})$.
(Caution: This is slightly different from the definition of $c_W$ 
in \cite{FSZ}.)

\subsection{Invariance constraints and universal relations} \label{s:5.3}
The discussion here in intended to get across general ideas, and will
be imprecise.  Interested readers may consult \cite{ypL4, LP} or
Givental's original papers (\cite{aG2} and references therein) for
details.

The starting point of this is Givental's axiomatic (semisimple) GWTs.
A major breakthrough in Givental's discovery is that the ``moduli space''
of Frobenius manifolds (or equivalently genus zero axiomatic theories) 
of \emph{a fixed dimension} has an action by a ``twisted loop group'',
or \emph{Givental's group}.
Furthermore, all semisimple theories form a single orbit under this group 
action.
One of the simplest semisimple theories of dimension $N$ is the GWT of 
$N$ points, which is given by $N$ copies of GWT of one point.
This implies that all semisimple axiomatic theories can be obtained
by a group action on GWT of $N$ points.
That is, for each semisimple axiomatic theory $T$, one can find an element
$O_T$ in Givental's group such that $T$ is equal to $O_T$ acting on 
the moduli point defined by GWT of $N$ points.
$O_T$ is uniquely determined if a ``homogeneity'' condition holds.
We shall not go further into the details, but only point out that
Witten's generating function $\tau_W$ satisfies the homogeneity condition.

A geometric GWT involves all genera and is encoded in a single generating 
function $\tau_{GW}$ as in \eqref{e:tauGW}.
For axiomatic theories, the generating functions $\tau_G$ are obtained via 
``quantizing'' $O_T$.
$\hat{O}_T$ is an operator which then acts on the generating function of
GWT of $N$ points $\tau^{N pt}$ as in  \eqref{e:tauG}.

One may define a \emph{universal relation} to be an equation of tautological
classes which holds for all axiomatic GWTs.
In particular it must contain tautological relations for moduli of curves as
it has to hold for $\tau^{N pt}$.
To show that a tautological relation for moduli curves holds for all axiomatic
GWTs, we can check whether it is ``invariant'' under the action of 
$\hat{O}_T$ for all axiomatic GWTs $T$.
In fact, since $O_T$'s form Givental group, 
we only have to check the invariance at the level of Lie algebra.
This gives very strong constraints on the form of the possible universal 
relations.

It was checked that all tautological equations are universal, \cite{ypL2}.
A simple geometric proof was later discovered by Faber--Shadrin--Zvonkine and
independently R.~Pandharipande and the first author 
(see \cite[\S 3]{FSZ}).

\subsection{The Faber--Shadrin--Zvonkine Uniqueness Theorem}

Let $\mu^i$ be partitions of $d$.
Let ${H}_{g}(\mu^1, \mu^2)$ be the moduli space of Hurwitz covers of degree 
$d$ from genus $g$ curves to $\mathbb{P^1}$, 
with fixed profiles $(\mu^1, \mu^2)$ at 
$(0, \infty)$ in $\mathbb{P^1}$, and otherwise simple ramification.
Let $\overline{H}_g(\mu^1, \mu^2)$ be the compactification by 
admissible covers.
There is a natural morphism
\[
 \rho: \overline{H}_g(\mu^1, \mu^2) \to \Mbar_{g, \sum_i l(\mu^i)}
\]
by forgetting the covering maps and stabilizing the domain curves.

One can consider variations of the above constructions.
First, one can allow partitions $\mu^i =(\mu^i_1, \mu^i_2, \ldots)$
to contain zeros, which correspond to marked unramified points.
Second, $\rho$ can be composed with forgetful maps 
(by forgetting marked points.)
For the purpose of this section, we are only interested in forgetting
points which are \emph{not} marked unramified points.
By abusing notation, we will denote them by
$\rho: \overline{H}_g \to \Mbar_{g,n}$.
Let $p := \sum_i l (\mu^i) -n$. 
(Note that a marked unramified point increases the length by 1.)
We will call the images of $\rho$  \emph{double Hurwitz cycles}.

Witten's correlators, by definition,
can be expressed as integrals of the following type
\[
  \int_{\Mbar_{g,n}} c_W \Psi,
\]
where $\Psi$ is a polynomial of $\psi$-classes.
Also, by construction, (complex) degree of $\Psi$ is at least $g$.
Now Ionel's theorem \cite{eI} states that any monomials in 
$\psi$ and $\kappa$ of Chow degrees greater than $g-1$ can be represented 
as a linear combination of classes of the form
\[
 q_*(DHC_1 \times DHC_2 \times \ldots)
\]
where $q: \prod_i \Mbar_{g_i, n_i} \to \Mbar_{g,n}$ is the 
gluing morphism and $DHC_i$ is a certain double Hurwitz cycle in 
$\Mbar_{g_i, n_i}$.
Since it is known \cite{FP} that double Hurwitz cycles are tautological,
and tautological classes are linearly generated by monomials of $\psi$ and 
$\kappa$ classes on (the vertices of) the dual graphs,
the above arguments imply that 
\emph{any monomial in $\psi$ and $\kappa$ on $\Mbar_{g,n}$ of 
Chow degrees at least $g$ can be represented as a linear combination of dual
graphs with at least one edge}.
This is called \emph{$g$-reduction} in \cite{FSZ}.

By $g$-reduction on $\Psi$ and the splitting principle of $c_W$,
Witten's correlators can then be reduced, with simple dimension counts and
manipulations of tautological classes, to genus zero correlators.

Note that the above arguments apply equally to axiomatic and geometric 
theories.
For axiomatic theories, one has to show that all tautological relations
holds, so that the same form of $g$-reduction applies.
As explained in the previous subsection, this is accomplished by showing
they satisfy the invariance constraints.

Now, the fact that both theories have 
\begin{enumerate}
 \item the same reduction to genus zero correlators, and 
 \item the same genus zero correlators
\end{enumerate}
implies the geometric theory equals the axiomatic one.
This, combined with Givental's theorem \cite{aG2} then proves
Witten's conjecture.

\end{document}